\documentclass[a4paper,12pt]{article}
\usepackage{inputenc}
\usepackage{graphics}
\usepackage{epsfig}
\usepackage{graphicx}
\usepackage{latexsym}
\usepackage{graphics}
\usepackage{epsfig}
\usepackage{amsmath,amssymb,amsfonts,theorem,color,bm}
\usepackage{multirow}
\usepackage{verbdef}
\usepackage{mathrsfs}
\usepackage{amsmath}
\usepackage{amsfonts}
\usepackage{epsfig}
\usepackage{booktabs}
\usepackage{graphics}
\usepackage{color}
\usepackage{ifthen}
\newcommand{\mA}{\mathsf{A}}
\newcommand{\mE}{\mathsf{E}}

\newboolean{showcomments}
\setboolean{showcomments}{false}

\usepackage{comment}
\usepackage{float}
\usepackage{subcaption}
\newcommand{\beqn}{\begin{equation}}
\newcommand{\eeqn}{\end{equation}}

\def\blackbox{\leavevmode\vrule height 5pt width 4pt depth 0pt\relax}

		\title{Modeling of air-wall heat transfer in buildings via recursive POD}

		\newcommand\restr[2]{{
				\left.\kern-\nulldelimiterspace 
				#1 
				\right|_{#2} 
			}}

\begin{document}
				\title{Anisotropic weights for RBF-PU interpolation with subdomains of variable shapes}
				
				\author{
				}
				\author{
					R. Cavoretto\thanks{Department of Mathematics \lq\lq Giuseppe Peano\rq\rq, University of Torino, Italy}
					\quad A. De Rossi\thanks{Department of Mathematics \lq\lq Giuseppe Peano\rq\rq, University of Torino, Italy}
					\quad G.~E. Fasshauer\thanks{Department of Applied Mathematics and Statistics,  Colorado School of Mines, Golden, CO, USA}
					\quad M.~J. McCourt\thanks{SigOpt, Inc., San Francisco, CA, USA} \\
					 E. Perracchione\thanks{Dipartimento di Matematica \lq\lq Tullio Levi-Civita\rq\rq, Universit\`a di Padova, Italy}
				}		
				
				\maketitle
				
				\section*{Abstract}
	The partition of unity (PU) method, performed with local
	radial basis function (RBF) approximants, has already been proved to be an effective tool for solving interpolation or collocation problems when large data sets are considered. It decomposes the original domain into several \emph{subdomains} or \emph{patches} so that only linear systems of relatively \emph{small} size need to be solved.
	In research on such partition of unity methods, such subdomains usually consist of spherical patches of a fixed radius. However, for particular data sets, such as track data,  ellipsoidal patches seem to be more suitable.
	Therefore, in this paper, we propose a scheme based on \emph{a priori} error estimates for selecting the sizes of such variable ellipsoidal subdomains.
	We jointly solve for both these domain decomposition parameters and the anisotropic RBF shape parameters on each subdomain to achieve
	superior accuracy in comparison to the standard partition of unity method.

\section{Introduction}

Radial basis function (RBF)-based methods \cite{Fasshauer15} find their natural applications in various fields, such as image reconstruction, resolution of partial differential equations and population dynamics.
Two common computational issues arising when we deal with real situations involve creating approximations from a
\emph{large} number of scattered points and the one of producing accurate approximations despite ill-conditioned linear systems.

In this article, we attack the first issue with an
efficient computation by means of the Partition of Unity (PU) method \cite{Wendland02a}. 
It enables us to decompose the original interpolation problem (involving a matrix of the same size as the amount of data) into many
small ones defined on subdomains/patches of the original domain.
However, the design of these subdomains (and consequently the number of points lying on each patch)
affects the accuracy of the approximation.

 {Generally, when the PU method is used for scattered data interpolation,
the PU subdomains are assumed to be balls of a fixed radius \cite{Cavoretto16}}.
\ifthenelse{\boolean{showcomments}}{
	{\color{red}[MM] Maybe we should add other PU references here, if we have space?}
}{}
In \cite{Cavoretto17}, a local approach was proposed via the PU method that
selects \emph{optimal local approximants}: both the shape parameter
and the patch radius were selected such that error estimates were minimized.
This strategy was more effective at dealing with points that were inconsistently distributed throughout the domain.

This previous work allowing varying patch radii was limited to only spherical patches.
Such a scheme is not completely suitable for particular data distributions, such as track data \cite{allasia}.
To address this situation we propose to use ellipsoidal subdomains, allowing one free domain parameter per dimension.
This adaptation requires careful selection of the PU weights to guarantee consistency, which we detail in our proposal.

To match the anisotropic structure of these subdomains, we use anisotropic Wendland's functions \cite{Wendland05} to
form our local approximants.
The values of the shape parameters and of the semi-axes of patches are selected by minimizing theoretical error estimates.
In particular, as in \cite{Cavoretto17}, we focus our attention on the Leave One Out Cross Validation (LOOCV) scheme \cite{Rippa}.
By jointly optimizing for both the domain and RBF parameters, we are able to improve on the computational cost from \cite{Cavoretto17}.
We also factor in the impact of ill-conditioning during this optimal parameter search to balance accuracy and stability.
\ifthenelse{\boolean{showcomments}}{
	{\color{red}[MM] We are saying that there is a speedup in computational time but we do not really see that in the results.  Do we need
		to add that comparison?}
	}{}

The outline of the paper is as follows. In Section \ref{remarks}, we briefly review the main theoretical aspects of the RBF-PU method. Section \ref{optimal} is devoted to the presentation of the proposed scheme which makes use of ellipsoidal patches. Numerical experiments are presented in Section \ref{numerical}. Section \ref{conclusions} deals with conclusions and work in progress.

\section{The RBF-based Partition of Unity method}
\label{remarks}

The approximation problem considered in this paper is formulated as follows.
Consider a set ${\cal X}_N = \{  \boldsymbol{x}_i,\; i = 1,\,\ldots,\,N \} \subseteq \Omega$  of distinct data points (or data sites or nodes), arbitrarily distributed on a domain $ \Omega \subseteq \mathbb{R}^{M}$, with an associated set $ {\cal F}_N= \{ f_i = f(\boldsymbol{x}_i) ,\;i=1,\,\ldots,\,N \}$ of data values (or measurements  or function values), which are obtained by sampling   some (unknown) function  $f: \Omega \longrightarrow \mathbb{R}$ at the nodes $ \boldsymbol{x}_i$.
The \emph{scattered data interpolation} problem consists of finding a function $R: \Omega \longrightarrow \mathbb{R}$ such that $R\left( \boldsymbol{x}_i\right)=f_i$, $i=1,\,\ldots,\,N$.

To this end, we take ${R} \in H_{\Phi} ({\cal X}_N)= \textrm{span} \{ \Phi(\cdot,\boldsymbol{x}_i),\; \boldsymbol{x}_i \in {\cal X}_N\}$, where $\Phi : \Omega \times \Omega \longrightarrow \mathbb{R}$ is a strictly positive definite and symmetric kernel.
More specifically, we take RBFs (radial kernels), and thus, we suppose that there exist a function $ \phi: [0, \infty) \to \mathbb{R}$ and a shape parameter $\varepsilon>0$ such that for all $\boldsymbol{x},\boldsymbol{y} \in \Omega$ we have $\Phi(\boldsymbol{x},\boldsymbol{y})=\phi_{\varepsilon}( ||\boldsymbol{x}-\boldsymbol{y}||_2):=\phi(r)$.
In Table \ref{tab_rbf}, we list the strictly positive definite RBFs that will be used later.
Note that the RBFs depend on a shape parameter $\varepsilon >0$ that significantly affects the accuracy of the approximation.
We will later refer to the functions reported in Table \ref{tab_rbf} as isotropic kernels, meaning that $\varepsilon$ is a scalar.

\begin{table}
	\caption{Examples of strictly positive definite isotropic radial kernels.}
	\label{tab_rbf}
	\begin{center}
	\begin{tabular}{cc}
		\hline\noalign{\smallskip}
		RBF  & $\phi(r)$ \\
		\hline\noalign{\smallskip}
		{\rm Inverse MultiQuadric $C^{\infty}$} (IMQ) & $(1+\varepsilon^2r^2)^{-1/2}$    \\
		{\rm Mat$\acute{\text{e}}$rn $C^2$} (M2) & ${\rm e}^{-\varepsilon r}    (\varepsilon r+1)$   \\
		{\rm Wendland $C^2$} (W2) & $\max \left(1-\varepsilon r,0\right)^4\left(4\varepsilon r+1\right)$  \\
		\noalign{\smallskip}\hline\noalign{\smallskip}
	\end{tabular}
		\end{center}
\end{table}

By using RBFs, the interpolant  assumes the form
\begin{equation}\label{eq0}
{R}(\boldsymbol{x}) = \sum_{k = 1}^N \alpha_k \phi(||\boldsymbol{x}- \boldsymbol{x}_k||_2), \quad \boldsymbol{x}\in\Omega.
\end{equation}
The coefficients $\boldsymbol{\alpha}= (\alpha_1,\,\ldots,\,\alpha_N)^T$ in \eqref{eq0} are determined by solving the linear system
$\mA \boldsymbol{\alpha}= \boldsymbol{f}$, where the entries of the matrix $\mA \in \mathbb{R}^{N \times N}$ are given by
$(\mA)_{ik}= \phi (||\boldsymbol{x}_i - \boldsymbol{x}_k||_2)$, $ i,k=1,\,\ldots,\,N$, and  $\boldsymbol{f}= (f_1,\,\ldots,\,f_N)^T$.
The uniqueness of the solution is ensured by the fact that the kernel $\Phi$ is strictly positive definite and symmetric.
\ifthenelse{\boolean{showcomments}}{
	{\color{red}[MM] I am not a fan of vectors with commas in them.  I don't really care, but I would write this as
		$\begin{pmatrix}\alpha_1&\cdots&\alpha_N\end{pmatrix}$.}
}{}

One drawback of this method is the computational cost associated with the solution of potentially large linear systems.
The PU method, presented below, enables us to overcome such issue.
At first, we cover the domain $\Omega$ with  $d$ overlapping {subdomains} $ \Omega_j$.   To be more precise, we require a regular covering, i.e., $ \{ \Omega_j \}_{j=1}^{d}$ must fulfill the following properties:
\begin{itemize}
	\item[i.] for each $ \boldsymbol{ x} \in \Omega$, the number of subdomains $ \Omega_j$, with $ \boldsymbol{x} \in \Omega_j$, is  bounded by a global constant $C_1$,
	\item[ii.] each subdomain $ \Omega_j$ satisfies an interior cone condition,
	\item[iii.] the local fill distances $ h_{ {\cal X}_{N_j}}$ are uniformly bounded by the global fill distance $h_{{\cal X}_N}$, where ${\cal X}_{N_j}= {\cal X}_N \cap \Omega_j$.
\end{itemize}

Once we select weight functions $W_j$, $j=1,\,\ldots,\,d$, the PU interpolant can be defined as
\begin{equation*}
{\cal I}\left( \boldsymbol{x}\right)= \sum_{j=1}^{d} R_j\left( \boldsymbol{x} \right) W_j \left( \boldsymbol{x}\right), \quad \textrm{with} \quad R_j\left( \boldsymbol{x}\right)= \sum_{k=1}^{N_j} \alpha_k^j \phi (||\boldsymbol{x} -  \boldsymbol{x}^j_k||_2),
\end{equation*}
where $R_j$ is defined on the subdomain $\Omega_j$, $N_j$ indicates the number of points on $\Omega_j$ and $\boldsymbol{x}^j_k \in {\cal X}_{N_j}$, with $k=1,\,\ldots,\,N_j$.
Therefore, the problem leads to solving $d$ linear systems of the form
$\mA_j \boldsymbol{\alpha}_j= \boldsymbol{f}_j$, where $\boldsymbol{\alpha}_j= (\alpha_1^j,\,\ldots,\,\alpha_{N_j}^j )^T$, $ \boldsymbol{f}_j = (f_1^j, \ldots , f_{N_j}^j )^T$ and the entries of $\mA_j \in  \mathbb{R}^{N_j \times N_j} $ are given by   $(\mA_j)_{ik}= \phi (||\boldsymbol{x}^j_i - \boldsymbol{x}^j_k||_2)$, $i,k=1,\,\ldots,\,N_j$.

Since the coefficients of the local interpolants are determined by imposing the local interpolation conditions, the functions $W_j$, $j=1,\ldots,d$, must form a partition of unity. Moreover, we also require that such partition of unity is  $k$-stable  {(see, e.g. \cite{Wendland05} Def.~15.16, p.~276)}, which in particular implies that $\textrm{supp}(W_j) \subseteq \Omega_j$. For instance, such conditions are satisfied for the well-known \emph{Shepard's weights}; refer e.g. to \cite{Wendland02a} for further details.
\ifthenelse{\boolean{showcomments}}{
	\textcolor{red}{[GF](As I understand Wendland's $k$-stable definition, this condition primarily requires boundedness of the derivatives of the weight functions)}
}{}
\ifthenelse{\boolean{showcomments}}{
	\textcolor{red}{[MM] I would probably cut this technical detail as it does not really provide much of a benefit here.}
}{}

\section{Optimal local interpolants for the RBF-based PU method}
\label{optimal}

We now focus on the selection of the PU patches, and we remove the standard assumption that they consist of balls of a fixed radius. Therefore we consider ellipsoidal patches, i.e., each $\Omega_j$ is defined through its semi-axes $\boldsymbol{\delta}_j=(\delta_1^j,\,\ldots,\,\delta_M^j)$. Moreover, in what follows,  we  use anisotropic kernels. We remark that any isotropic radial kernel can be turned into an anisotropic one by using a weighted $2$-norm instead of an unweighted one. Thus, to fix the ideas, on a subdomain $\Omega_j$ it is enough to replace the scalar value of the shape parameter $\varepsilon_j$ with a symmetric positive definite matrix $\mE_j$. More precisely, we  consider the special case for which $\mE_j = \textrm{diag} (\varepsilon^j_1,\,\ldots,\,\varepsilon^j_M)$. This allows us to choose a different scaling along the dimensions of the problem.

Our parametrization strategy consists of selecting both  $\boldsymbol{\varepsilon}_j=(\varepsilon^j_1,\,\ldots,\,\varepsilon^j_M)$ and $\boldsymbol{\delta}_j$ such that the error estimates on $\Omega_j$ are minimized. This study is motivated by the fact that the PU approximation error is governed by the local ones (\cite{Wendland05} Th.~15.19, p.~277).

\subsection{Local error estimates}

For a general overview about error estimates refer e.g. to \cite{Fasshauer15}. Here we focus on error predictions that are popular in statistics, and precisely on cross validation schemes. We describe the cross validation algorithm that is applied on a given $\Omega_j$  to get a local \emph{a priori} error estimate \cite{Fasshauer15}. At first, we split the set ${\cal X}_{N_j}$ into two disjoint subsets: a training set ${\cal T}_{N^t_j}$ and a validation set  ${\cal V}_{N^v_j}$ such that $N^t_j+N^v_j={N}_j$. The set  ${\cal T}_{N^t_j}$  is used to construct a surrogate or partial approximation that is validated via the set ${\cal V}_{N^v_j}$. To simplify the following discussion, on $\Omega_j$ we introduce the following block decomposition of the local interpolation matrix
\begin{align*}
\mA_j = \left(
\begin{array}{ccc}
\mA_j^{tt}  & \mA_j^{tv}       \\
\mA_j^{vt}  & \mA_j^{vv}
\end{array}
\right),
\end{align*}
where, for example, the block $\mA_j^{tv}$ is generated using training points to evaluate and validation data as centers of the kernels. Similarly, we partition ${\boldsymbol{\alpha}}_j = (\boldsymbol{\alpha}_j^{t}, \boldsymbol{\alpha}_j^{v})^T$, $\boldsymbol{f}_j = (\boldsymbol{f}_j^{t},
\boldsymbol{f}_j^{v} )^T$.
With this notation, the prediction at the points in the validation set using the training set is $\mA_j^{vt} (\mA_j^{tt})^{-1} \boldsymbol{f}_j^{t}$.
In other words, $|\boldsymbol{f}_j^{v}-\mA_j^{vt} (\mA_j^{tt})^{-1} \boldsymbol{f}_j^{t}|$ provides information about the accuracy of the fit on the  $j$-th subdomain.
\ifthenelse{\boolean{showcomments}}{
	\textcolor{red}{(GF comment: I've changed the notation for the following discussion. However, this discussion with $q$ partitions is more general than leave ONE out. I'm not sure if we really should be putting this here. Probably we should just keep things at LOOCV.)}
}{}
Then, we take $q$ partitions of  ${\cal V}_{N^v_j}$, ${\cal V}_{N^v_j} = \{{\cal V}_{N^v_j}^{(1)}, \ldots , {\cal V}_{N^v_j}^{(q)} \}$, such that
\begin{equation*}
	{\cal V}_{N^v_j}^{(k)} \cap {\cal V}_{N^v_j}^{(i)} = \emptyset, \hskip 0.2cm
	\textrm{for} \hskip 0.2cm i \neq k, \hskip 0.2cm
	\textrm{with} \hskip 0.2cm \displaystyle{\cup_{k=1}^{q} }{\cal V}_{N^v_j}^{(k)} = {\cal X}_{{{N}_j}}, \hskip 0.2cm
	\textrm{and} \hskip 0.2cm {\cal T}_{N^t_j}^{(k)}={\cal X}_{{{N}_j}}^{(k)} \backslash {\cal V}_{N^v_j}^{(k)}.
\end{equation*}
Thus, as error estimate we  consider the residual left over by
the interpolants evaluated at the validation sets.  {The LOOCV scheme is a particular case of the general setting presented above for which $q=N_j$ and each ${\cal V}_{{N^v_j}}^{(k)} = \boldsymbol{x}_k^j$.}
\ifthenelse{\boolean{showcomments}}{
	\textcolor{red}{(GF comment: If we do want to keep the general discussion with the $q$ partitions, then writing something like ``The LOOCV scheme is a particular case of the general setting presented above for which $q=N_j$ and each ${\cal V}_{\textcolor{red}{N^v_j}}^{(k)} = \boldsymbol{x}_k^j$'' would probably be better.}
}{}
Moreover  for the LOOCV scheme, $\mA_j^{vv}$  is the diagonal element of $\mA_j^{-1}$ and thus, being a scalar, the computation simplifies. Indeed, as error estimate for the $j$-th subdomain we have (see also \cite{Rippa})
\begin{equation*}
e_j=  ||(\alpha^j_1 / (\mA_j^{-1})_{11},\,\ldots,\,\alpha^j_{N_j } / (\mA_j^{-1} )_{N_j N_j}) ||_p,
\end{equation*}
where in what follows we fix the index of the discrete norm to $p=2$.

In our PU context, using both anisotropic kernels and ellipsoidal patches, we have that $e_j=e_j(\boldsymbol{\varepsilon_j},\boldsymbol{\delta_j})$.
In fact, the shape parameter affects the accuracy of the RBF approximant and, for the PU method, the accuracy also depends on which points are involved in the computation of the local interpolants.

\subsection{Description of the PU-LOOCV method}

To minimize the LOOCV error estimates we use a multivariate optimization tool. This allows us to reduce the computational cost of the procedure presented in \cite{Cavoretto17}.  To be more precise, we consider  the Nelder-Mead simplex algorithm \cite{Lagarias}. Without going into details, we remark that on $\Omega_j$, given an initial guess $(\boldsymbol{\varepsilon}_j^0, \boldsymbol{\delta}_j^0)$, it provides an approximation of the optimal values by computing subsequent simplices and  only needs function evaluations of the objective function \cite{Lagarias}. In particular,  for the implementation we use the \textsc{Matlab} software and the  {\tt fminsearch.m} routine.
Of course, we also need to impose the following constraints on the parameters we optimize:
\begin{equation*}
{\varepsilon}_k^j > 0, \quad \textrm{and} \quad  \delta^+ \geq {\delta}_k^j \geq \delta^*, \quad k=1,\ldots,M, \quad j=1,\ldots,d,
\end{equation*}
where $\delta^* \in \mathbb{R}$ is chosen so that patches form a covering of the domain and $\delta^+ \in \mathbb{R}$ is selected so that for each $ \boldsymbol{ x} \in \Omega$ the number of subdomains $ \Omega_j$, with $ \boldsymbol{x} \in \Omega_j$, is bounded.

In this way, after optimizing the parameters $(\boldsymbol{\varepsilon}_j, \boldsymbol{\delta}_j)$, $j=1,\ldots,d$, we have a PU covering made of ellipsoids that is also regular. In fact, each of the patches satisfies an interior cone condition. This assumption is trivially verified for balls. However, it is also true for ellipsoids (see, e.g. \cite{Wendland05} Pr.~11.26, p.~195).

Finally, in order to build a consistent PU setting, we also need to carefully choose the compactly supported functions for the PU Shepard's weights, which are here constructed  with the W2 function. Since we have ellipsoids, we need to select the anisotropic form of the W2 function. Thus, $\boldsymbol{\varepsilon}_j$, which identifies the support of the compactly supported RBF, is taken so that ${ \rm supp}(W_j)  = \Omega_j$,  $j=1, \ldots,d$.

\section{Numerical experiments}
\label{numerical}

Our experiments focus on bivariate interpolation. In Subsection \ref{ExArtificial}, we show the numerical results obtained by considering known functions and artificial track data \cite{allasia}. Then, in Subsection \ref{ExReal} we also take into account real data by analyzing an application to Earth's topography.

\subsection{Experiments with artificial data}
\label{ExArtificial}

To illustrate the accuracy of the proposed method, we evaluate the interpolant on a grid of $s=40^2$  points $\tilde{\boldsymbol{x}}_i$, $i=1,\,\ldots,\,s$, on $\Omega=[0,1]^2$ and we calculate Root Mean Square Error (RMSE).
The patch centers are constructed as a grid of $t^2$ points on $\Omega$, where $t$ is the number of tracks. Of course this design will affect the accuracy of the approximation. Nevertheless, the scheme proposed here, allowing to choose variable subdomains, is consequently less sensitive to this choice.
\ifthenelse{\boolean{showcomments}}{
	\textcolor{red}{(GF comment: I guess we're not going to say more about this here, but this choice of (distribution of) patch centers deserves further investigation.)}
}{}

We show  numerical results obtained by considering five sets of track data on $\Omega=[0,1]^2$ sampled {from} the \emph{2D Franke's function}, see e.g. \cite{Fasshauer15}. In particular, the results of using LOOCV to optimize the semi-axes of the patches and the shape parameters of the local basis functions are reported in Table \ref{tab_1}.  {They are compared with the classical PU method obtained by taking a grid of $t^2$ points on $\Omega$ as PU centers and a fixed patch radius $\delta=\delta^*$. Furthermore, we consider the isotropic IMQ kernel with shape parameter equal to $1$. We select such shape parameter arbitrarily. Indeed, one of the main advantages of the proposed method is the one of automatically choosing \emph{safe} values for the shape parameters. Of course, different values might lead to more accurate approximations, but it is not possible to  provide a priori optimal or safe shape parameters.} Moreover, we also report the CPU times. Tests have been carried out with the \textsc{Matlab} software on a{n} Intel(R) Core(TM) i7-6500U CPU 2.59 GHz processor. In Figure \ref{fig_1} (left), we show an example of $2000$  track data \cite{allasia} and the ellipsoidal patches obtained via PU-LOOCV.

From the numerical experiments we can note that the classical PU method, which makes use of circular patches, is not able to accurately fit the data, especially when a large number of points is involved.  {This might be due to a non-optimal selection of the shape parameter and/or of the patch size for the classical PU. In this sense, the PU-LOOCV reveals its robustness, selecting optimal values for those parameters and providing  accurate results also when $N$ grows.} Finally, we point out that the proposed scheme, besides extending the idea of the method presented in \cite{Cavoretto17} to subdomains having different shapes and to anisotropic kernels, thanks to the use of an optimization routine for the minimization problem, it also speeds up the procedure.
For instance, using the same scheme outlined in \cite{Cavoretto17} would take about $300$ s for $1000$ data. Moreover, note that Table \ref{tab_1} shows that the PU-LOOCV with $1000$ points is more than twice as accurate \emph{and} more than twice as fast than the classical PU method with $160000$ points.

\ifthenelse{\boolean{showcomments}}{
	{\color{red}[MM] Is there any benefit to providing both the RMSE and the MAE?  I do not really see that there is.  I think the table
		might be more readable as I have it below.  But that is just my opinion.}
}{}

\begin{table}[ht]
	\centering
	\caption{RMSEs, and CPU times obtained by using the PU-LOOCV and classical PU methods.\label{tab_1}}
	\setlength{\tabcolsep}{10pt}
	\begin{tabular}{lrrrr}
		\hline\noalign{\smallskip}
		\multicolumn{1}{c}{$N$}  & \multicolumn{2}{c}{RMSE} & \multicolumn{2}{c}{CPU time}  \\ \cmidrule(lr){2-3}\cmidrule(lr){4-5}
		& PU & PU-LOOCV & PU & PU-LOOCV  \\
		\noalign{\smallskip}\hline\noalign{\smallskip}
		$1000$ $(20\times50)$   &  $3.27{\rm E}\mbox{--}3$  & $8.02{\rm E}\mbox{--}5$  &  $0.5$  & $43.3$  \\
		$2000$ $(25\times80)$   &  $1.08{\rm E}\mbox{--}3$  & $2.58{\rm E}\mbox{--}5$  &   $1.2$  & $64.7$  \\
		$4000$ $(40\times100)$   &  $1.28{\rm E}\mbox{--}3$  & $4.67 {\rm E}\mbox{--}6$  &   $7.2$  & $136.0$  \\
		$8000$ $(50\times160)$   &  $3.48{\rm E}\mbox{--}4$  & $1.25 {\rm E}\mbox{--}6$  &   $17.4$  & $189.0$  \\
		$16000$ $(80\times200)$   &  $3.49{\rm E}\mbox{--}4$  & $4.55{\rm E}\mbox{--}7$  &   $109.0$  & $475.0$  \\
		\noalign{\smallskip}\hline\noalign{\smallskip}
	\end{tabular}
\end{table}

\subsection{Experiments with real data}
\label{ExReal}

\begin{figure}[t]
	\centering
	\includegraphics[width=.46\linewidth]{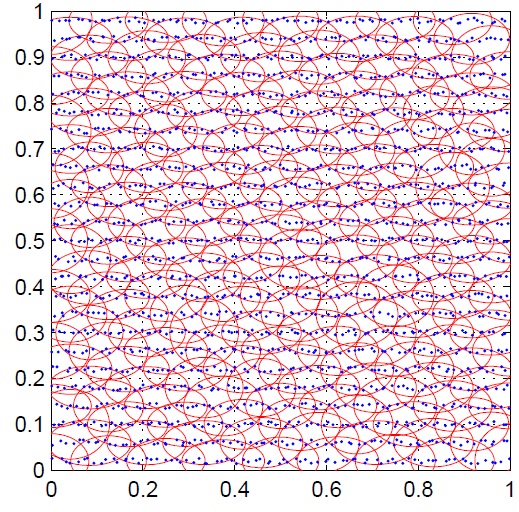}
	\hspace{.04\linewidth}
	\includegraphics[width=.43\linewidth]{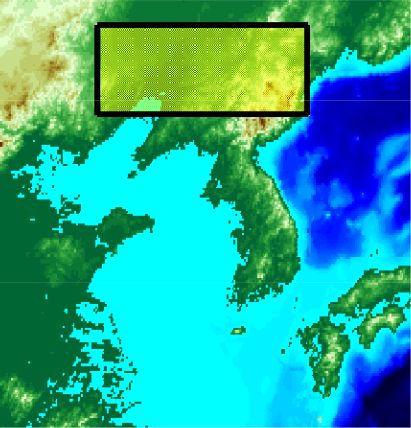}
	\caption{Left: an illustrative example with $2000$ track data that shows how patches are selected by means of the PU-LOOCV method.
		Right: the Korea's map and the extracted tracks.
		\label{fig_1}
	}
\end{figure}

To test the method with real data, we consider points extracted from maps. We take, as example, {a map of Korea} (plotted in the right frame of Figure \ref{fig_1}) and we extract $40$ tracks containing $100$ points.   {The function values, being real samples of the elevation above sea level, are truly oscillating  and thus the interpolation problem is particularly challenging.} We also extract a grid of $40^2$ points to evaluate the error.

 {By using the M2 kernel,  the RMSE for the PU-LOOCV is equal to $2.58{\rm E}-2$. The classical PU with a fixed radius $\delta=\delta^*$ completely fails. However, for the classical PU, taking anisotropic kernels with variable parameters and $\delta= 3 \delta^*$ as fixed size of the circular patches gives RMSE=$4.21{\rm E}-2$. In other words, we reach about the same accuracy, but this is not completely satisfying. Indeed, if we take circular subdomains of radius $3 \delta^*$, the average of points on each patch is about $52$, while the PU-LOOCV only requires on average $16$ data per patch. Therefore, especially with real data, the approximation by means of an optimized PU method becomes essential.}

\section{Conclusions}
\label{conclusions}

In this paper we presented a scheme for the optimal selection of local approximants in the PU method. The scheme, based on ellipsoidal patches, is particularly suitable for track data. Work in progress consists in comparing the LOOCV scheme with other a priori error estimates  {and in selecting suitable locations for the patch centers. Indeed, here they are taken as grids of points, but this might be restrictive.} Further investigations for more efficient optimization routines are also needed. 

\section*{Acknowledgments}     	
This research has been accomplished within RITA (Rete ITaliana di Approssimazione) and partially supported by GNCS-IN$\delta$AM. The first and second authors were partially supported by the 2016-2017 project \emph{Metodi numerici e computazionali per le scienze applicate} of the Department of Mathematics of the University of Torino. The third author was partially supported by grant NSF-DMS \#1522687. The last author is supported by the research project \emph{Radial basis functions approximations: stability issues and applications}, No. BIRD167404.


\end{document}